\documentclass[a4paper]{article}
\usepackage{amsmath,amssymb, amsfonts, mathbbol, mathrsfs}

\usepackage{RR}
\RRNo{6255}

\usepackage{hyperref}

\RRdate{July 2007}

\RRauthor{ Bruno Cernuschi-Fr{\'\i}as\thanks{University of Buenos Aires, and CONICET, Argentina}}

\authorhead{Cernuschi-Fr{\'\i}as}

\RRtitle{Champs de Markov al{\'e}atoires {\`a} {\'e}tats mixtes avec {\'e}tiquettes symboliques et valeurs multi-dimensionnelles r{\'e}elles}

\RRetitle{Mixed States Markov Random Fields with Symbolic Labels and Multidimensional Real Values}

\titlehead{Mixed States MRF's}

\RRnote{Joint collaboration between the  Equipe  Associ{\'e} INRIA {\bf FIM}, {\it Fluides, Images et Mouvement}, CONICET, and Facultad de Ingenier{\'\i}a, Universidad de Buenos Aires, Argentina, 
and VISTA, IRISA - INRIA, Rennes.}

\RRresume{
De nouveaux r{\'e}sultats th{\'e}oriques sont present{\'e}s ici sur le mod{\`e}le r{\'e}cemment introduit appel{\'e} Champs de Markov Al{\'e}atoires 
{\`a} {\'E}tats Mixtes. 
Ces mod{\`e}les ont {\'e}t{\'e} introduits dans le contexte de l'analyse du mouvement dans des s{\'e}quences d'images, et sont utiles pour repr{\'e}senter l'information qui peut prendre des valeurs discr{\`e}tes correspondant {\`a} des {\'e}tats symboliques, et des valeurs r{\'e}elles correspondant {\`a}
des mesures continues. 
En particulier, des r{\'e}sultats sont donn{\'e}s quand l'{\'e}nergie globale pour la formulation de Gibbs, exprimant le mod{\`e}le {\`a} {\'e}tats mixtes, 
peut {\^e}tre d{\'e}compos{\'e}e en une partie discr{\`e}te et une partie continue. Cette d{\'e}composition permet de d{\'e}finir les mod{\`e}les conditionnels {\`a} {\'e}tats  mixtes d'une mani{\`e}re tr{\`e}s simple.

}

\RRabstract{
New theoretical results are presented here on the recently introduced
model called “mixed states MRF”. Such models were
introduced in the context of image motion analysis and are
useful to represent information which can take both discrete
values accounting for symbolic states, and real values
corresponding to continuous measurements. In particular, results are given when
the global energy for the Gibbs formulation expressing the mixed states model,
can be decomposed into one term accounting for the
discrete part of the model, and a second term related to the continuous
part. This decomposition theorem permits
to define conditional mixed states models in a very simple way. 
}

\RRmotcle{Analyse de mouvement en images, champs de Markov al{\'e}atoires, {\'e}tats mixtes.}

\RRkeyword{
Image motion analysis, Markov Random Fields, Mixed States.
}

\RRprojet{VISTA}

\RRtheme{\THCog}

\URRennes  

\newtheorem{theorem}{Theorem}[section]
\def\stackblw#1#2{\mathrel{\mathop{#1}\limits_{#2}}}

\newcommand{\SPO}{\Omega}
\newcommand{\SAF}{\mathcal{F}}
\newcommand{\REALS}{\mathbb{R}}
\newcommand{\BOREL}{\mathscr{B}}
\newcommand{\MUX}{\mu_X}
\newcommand{\MMUX}{\MUX^{m}}
\newcommand{\MDX}{m^{d}}
\newcommand{\DFX}{F_X}
\newcommand{\PDFX}{p_X^{a}}
\newcommand{\MPDFX}{p_X^{m}}
\newcommand{\STEP}{u}
\newcommand{\UNOG}{\mathbb{1}}
\newcommand{\PM}{\mathscr{P}}

\newcommand{\MRFX}{p_{\bf X}^{m}}

\newcommand{\SPL}{\mathbb{L}}

\newcommand{\ELE}{\mathbb{L}}
\newcommand{\EME}{\mathbb{M}}
\newcommand{\SAM}{\mathcal{M}}

\newcommand{\XRV}{{\bf X}}
\newcommand{\PDFXV}{p_{\XRV}^{a}}
\newcommand{\MPDFXV}{p_{\XRV}^{m}}
\newcommand{\MUXV}{\mu_{\XRV}}
\newcommand{\MMUXV}{\MUXV^{m}}

\newcommand{\DELTA}{\UNOG}
\newcommand{\LOG}{\log}
\newcommand{\PC}{p^{{a}}}

\newcommand{\FC}{Q^{{a}}}
\newcommand{\FD}{Q^{{d}}}
\newcommand{\FMS}{Q^{m}}
\newcommand{\VD}{V^{d}}
\newcommand{\VC}{V^{a}}

\newcommand{\VDX}{V^{d}(\X)}
\newcommand{\VCX}{V^{a}(\X)}
\newcommand{\VMX}{V^{m}(\X)}
\newcommand{\X}{{\bf x}}
\newcommand{\XA}{\X_{A}}
\newcommand{\R}{{\bf r}}
\newcommand{\FMSX}{Q^{m}(\X)}
\newcommand{\FCX}{\FC(\X)}
\newcommand{\FDX}{\FD(\X)}
\newcommand{\PCX}{\; \PC(\X) \;}

\newcommand{\PCR}{\; \PC(\R) \;}
\newcommand{\LFPC}{\; \LOG \frac{\PCX}{\PCR} \;}
\newcommand{\XK}{{\bf x}_{i}}
\newcommand{\XL}{{\bf x}_{l}}
\newcommand{\XLRV}{{\bf X}_{l}}
\newcommand{\XJ}{{\bf x}_{j}}
\newcommand{\XM}{{\bf x}_{k}}

\newcommand{\RK}{{\bf r}_{i}}
\newcommand{\RL}{{\bf r}_{l}}
\newcommand{\RKC}{{\bf r}_{i}^{c}}

\newcommand{\RJ}{{\bf r}_{j}}
\newcommand{\RM}{{\bf r}_{k}}
\newcommand{\XKC}{\XK^{c}}
\newcommand{\XLC}{\XL^{c}}
\newcommand{\XLCRV}{\XLRV^{c}}
\newcommand{\PMS}{p^{m}}
\newcommand{\PMSX}{\; \PMS(\X) \;}
\newcommand{\PMSR}{\; \PMS(\R) \;}
\newcommand{\LFPMS}{\; \LOG \frac{\PMSX}{\PMSR} \;}
\newcommand{\PKMS}{\; \PMS(\XK  \; \vert \; \XKC) \;}
\newcommand{\LPKMS}{\; \LOG \PKMS \;}
\newcommand{\QK}{\; \rho_{i} \left(\XKC \right) \;}
\newcommand{\QKS}{\; \rho^{*}_{i} \left(\XKC \right) \;}
\newcommand{\LQK}{\; \LOG \QK \;}
\newcommand{\LQKS}{\; \LOG \QKS \;}
\newcommand{\DK}{\; \DELTA_{\RK}(\XK) \;}
\newcommand{\DKS}{\; \DELTA_{\RK}^{*}(\XK) \;}
\newcommand{\DJS}{\; \DELTA_{\RJ}^{*}(\XJ) \;}
\newcommand{\DMS}{\; \DELTA_{\RM}^{*}(\XM) \;}
\newcommand{\PCK}{\; \PC \left(\XK  \; \vert  \; \XKC \right) \;}
\newcommand{\PCL}{\; \PC_{\XLRV \; \vert \; \XLCRV} \left(\XL  \; \vert  \; \XLC \right)}
\newcommand{\LPCK}{\; \LOG \PCK \;}
\newcommand{\PCRK}{\; \PC (\RK  \; \vert \;  \XKC) \;}
\newcommand{\PCRL}{\; \PC_{\XLRV \; \vert \; \XLCRV} (\RL  \; \vert \;  \XLC) \;}
\newcommand{\LPCRK}{\; \LOG \PCRK \;}

\newcommand{\LFCK}{\; \LOG \frac{\PCK}{\PCRK} \;}

\newcommand{\NLFCL}{\; \frac{\PCL}{\PCRL} \;}
\newcommand{\CCK}{\; \sum_{\{i\} \subseteq A \subseteq S} \FC_{A}(\XA) \;}
\newcommand{\CCL}{\; \sum_{\{l\} \subseteq A \subseteq S} \FC_{A}(\XA) \;}
\newcommand{\CC}{\; \sum_{A \subseteq S} \FC_{A}(\XA) \;}
\newcommand{\CMSK}{\; \sum_{\{i\} \subseteq A \subseteq S} \FMS_{A}(\XA) \;}
\newcommand{\CMS}{\; \sum_{A \subseteq S} \FMS_{A}(\XA) \;}
\newcommand{\CHC}{\; \sum_{A \subseteq S} Q_{A}(\XA) \;}
\newcommand{\CHCC}{\; \sum_{\{i\} \subseteq A \subseteq S} Q_{A}(\XA) \;}
\newcommand{\BGK}{ \mathfrak{G}^{(i)}}
\newcommand{\BGJ}{ \mathfrak{G}^{(j)}}
\newcommand{\BGM}{ \mathfrak{G}^{(k)}}
\newcommand{\CGK}{\; \sum_{ A \subseteq S \setminus \{i\}} \BGK_{A}(\XKC) \;}
\newcommand{\LQSP}{\; \LOG \left( \QKS \; \PCRK \right) \;}
\newcommand{\LFQP}{\; \LOG \left[ \frac{\QKS}{\QK}\PCRK \right] \;}
\newcommand{\PRKMS}{\; \PMS(\RK \; \vert \; \XKC) \;}
\newcommand{\LPRKMS}{\; \LOG \PRKMS \;}
\newcommand{\LFPKMS}{\; \LOG \frac{\PKMS}{\PRKMS} \;}
\newcommand{\GKXC}{\; \mathfrak{g}_{i}(\XKC) \;}
\newcommand{\GKR}{\; \mathfrak{g}_{i}(\RKC) \;}

\newcommand{\FOCKMS}{\; \FMS_{\{i\}}(\XK) \;}

\newcommand{\QKR}{\; \rho_{i}(\RKC) \;}

\newcommand{\QKRS}{\; \rho^{*}_{i}(\RKC) \;}

\newcommand{\PCRR}{\; \PC ( {\bf r}_{i} \; \vert \; {\bf r}_{i}^{c})}

\newcommand{\FCCK}{\; \FC_{\{i\}} (\XK)}
\newcommand{\AK}{\; \alpha_{i}}
\newcommand{\AKDEF}{\; \LOG \left[ \frac{\QKRS}{\QKR} \PCRR \right]}

\newcommand{\ALPHAD}{{\boldsymbol \phi}_{{d}}}
\newcommand{\ALPHAC}{{\boldsymbol \phi}_{{a}}}
\newcommand{\FCXP}{\VC(\X ; \ALPHAC)}
\newcommand{\PDXP}{p^{d}(\X ; \ALPHAD)}
\newcommand{\PCXP}{p^{a}(\X ; \ALPHAC)}
\newcommand{\FDXP}{\VD(\X ; \ALPHAD)}

\newcommand{\ZD}{Z^{{d}}(\ALPHAD)}
\newcommand{\ZC}{Z^{{a}}(\ALPHAC)}
\newcommand{\ZMS}{Z^{m}(\ALPHAD,\ALPHAC)}
\newcommand{\FMSXP}{V^{m}(\X;\ALPHAD,\ALPHAC)}
\newcommand{\ZN}{Z(\ALPHAD,\ALPHAC)}
\newcommand{\HKC}{\mathfrak{h}_{i}(\XKC)}
\newcommand{\HLC}{\mathfrak{h}_{l}(\XLC)}

\begin{document} 

\makeRR

\section{Introduction}\label{sec:intro}

Usually in statistics one is either interested in random variables that are discrete in nature, (dices, card games, photon counting, nuclear decay, etc), or
that are continuous, meaning that its distribution function is absolutely continuous, and hence they have a probability density function defined as the derivative of the distribution function, as is the case, for example, of the Gaussian random variables.
In some situations, one may be interested in modeling a variable that may take some discrete values in a set with non-zero probability mass, while for the other values not in the given set, the random variable may be modeled as a continuous distribution. Examples, might be the daily rainfall, 
\cite{Hardouin-1}, for which there are days with no rain with non zero probability,
while for the rainy days a continuous distribution random variable might be appropriate. In quantum mechanics one is faced with such situations, as for example in solid state-physics when one considers the solution of the band model for semiconductors. Another example arises in the field of reliability when modeling the mean life time of some component, for which there is a non zero mass probability of failure at time zero, while for time greater than zero an absolutely continuous distribution might be appropriate. In particular, these models appear to be of interest in low level modeling
of motion in image processing. In previous work, \cite{Bouthemy-1}, \cite{SPIE-06}, \cite{ICIP-06}, \cite{Gwenaelle-1}, \cite{Gwenaelle-2}, 
\cite{Gwenaelle-3}, \cite{Gwenaelle-4}, \cite{Gwenaelle-5}, \cite{yao-1}, it has been found that when modeling motion in complex images, as for example a video sequence showing public in a stadium, or a scene of a tree waving under the wind, among other examples, the histogram of the velocity vectors at each pixel show a large
probability mass at zero velocity, while the second component of the mixture may be appropriately modeled with a Gaussian distribution in many situations.

These examples suggest the introduction of mixed states random variables, \cite{Bouthemy-1}, \cite{SPIE-06}, \cite{ICIP-06}, \cite{Gwenaelle-5}. 
That is, variables
that have mass probability concentrated in isolated points, while they have a probability density for the rest of the real line.

Image motion analysis involves challenging issues such as
motion detection, segmentation, estimation, recognition or classification,
\cite{cedras},  \cite{radke}, \cite{stille}. In this context, compact and efficient representations
of image motion are needed.
Digital information present in or extracted from images
may be expressed as numerical values or discrete values (i.e.,
abstract labels). Moreover, these two types
of variables may more deeply reflect two different classes of information:
continuous real values (either in one-dimensional
or multi-dimensional spaces) versus symbolic values (one or
several symbols). However, these two classes should not be
necessarily viewed as two exclusive states or as consecutive states
(e.g., after a decision step). Indeed, a physical variable can
take both continuous and discrete values, namely it can be a
mixed states variable.
To give a simple example related to image motion analysis,
a locally computed motion quantity can be either null
or not. Then, it can be helpful to explicitly consider that it
takes either a discrete value expressing the absence of motion,
or that continuous real values account for the actual measurements.
A discrete value may not necessarily be a specific real value, it may also be taken as a pure symbolic value as well.
As an example, when considering optical flow, the label could be related to the
presence of motion discontinuities while the corresponding
continuous values are the velocity vectors. Evaluating the distribution
of measured values, and accounting for local context,
are of key importance in numerous image sequence analysis
tasks (e.g. in motion modeling, detection, segmentation, estimation,
recognition, or learning issues). Therefore, defining
probabilistic mixed states models appears as an attractive objective.
Markov Random Fields (MRF) are widely used in image
analysis, \cite{chalmond}. Recently, the so called mixed states
auto-models were introduced, \cite{Bouthemy-1}, (i.e., MRF models with two-site cliques), as a generalization of
the models studied in \cite{BESAG-1}.
These models were used for modeling and segmenting
motion textures with convincing experimental results, \cite{ICIP-06}. 

New theoretical
results, in the form of Gibbs potentials, allow to significantly extend
the power expression and the effective and flexible
use of these models, that may open new investigation avenues in image motion
analysis tasks. 

The contributions presented here are made along the following lines:

\begin{enumerate}

\item \label{it-1}
A detailed measure theoretic formulation of the mixed states random variables is given.

\item \label{it-2}
The cases in which the probability mass is concentrated 
either in discrete real values and/or 
in symbolic labels is thoroughly analyzed.

\item \label{it-3}
A full account is given for Markov Random Fields of mixed states random variables corresponding 
to a mixture of a probability mass in a known real value, and an absolutely continuous distributed 
real random variable as second component of the mixture, when there is 
no interaction between the potentials of the discrete and the absolutely continuous components of the Gibbs representation.

\item \label{it-4}

The generalized Markov Random Fields presented here are of general form in the sense that no stationarity assumption is assumed, 
moreover, at each of the sites the corresponding random variables may take values in different spaces.

\item \label{it-5}
The theoretical approach introduced here permits to analyze mixed states Markov Random Fields for any clique size, 
hence extending the previous results
given in 
\cite{Bouthemy-1}, \cite{Hardouin-1}, \cite{Hardouin-2}, \cite{Gwenaelle-5},
which are only valid for random variables belonging 
to the exponential family and up to cliques of size two.

\end{enumerate}

\section{Mixed States Random Variables} \label{MS-RVS}

\subsection{Preliminaries}\label{sec:prelims}

It is well known that the distribution function of a random variable
taking one-dimensional real values may be decomposed in the convex combination of three parts: 
a discrete part, an absolutely continuous part, and a singular continuous part, 
all with respect to the Lebesgue measure. 
The random variables considered here have as distribution function a convex combination of a discrete and an
absolutely continuous (with respect to Lebesgue measure) distribution functions. These distribution functions do not have a probability density 
function with respect to the Lebesgue measure, but they  are absolutely continuous with respect
to a measure which is the sum of the unit mass probability at the discrete points plus the 
Lebesgue measure over the real line. Call the corresponding Radon-Nikodym derivative a "mixed states probability density function
(ms-pdf)" and the corresponding random variable, a "mixed states random variable". This idea permits to immediately consider random variables on more general spaces. In particular, as previously discussed, it is of interest to study 
random variables that take values either in a countable set of labels $\SPL$ 
with probability $0 < \rho < 1$, or in $\REALS^{n}$ with probability $1 - \rho$. 
That is, given a probability space $(\SPO, \SAF, \PM)$, it is of interest to analyze 
random variables of the form $X: \SPO \rightarrow \SPL \cup \REALS^{n}$.
 
The following convention is used for upper scripts of probability related functions: the upper script $^{m}$ will correspond to mixed states related formulations,
the upper script $^{d} $ will correspond to discrete states formulations, while $^{a}$ will correspond to formulations related to absolute continuity with respect to the Lebesgue measure.

\subsection{Distribution function for one dimensional real random variables}\label{sec:DF-REAL}

Let's start with a brief review of distribution functions for one dimensional real random variables.

The following results may be found in any standard text on probability theory such as \cite{CHUNG}, 
as well as any standard text on measure theory such as \cite{HALMOS} or \cite{H-S}.

Consider a finite {\em measure space} $(\SPO,\SAF,\PM)$ where $\SPO$ is some abstract space, $\SAF$ is a $\sigma$-algebra of subsets 
of $\SPO$ and $\PM$ is a finite {\em measure} for the {\em measurable space} $(\SPO,\SAF)$. A {\em probability space} is a finite measure
space with $\PM(\SPO) = 1$.

Call $\mathbb{N}$ the set of non-negative integers.
Call $\REALS$ the set of the real numbers, and call  $\BOREL$ the Borel $\sigma$-algebra of subsets of $\REALS$, 
i. e. the minimal $\sigma$-algebra generated, say, by the open subsets of $\REALS$.

A {\em real finite-valued random variable} $X$ is a function from $\SPO$ to the real numbers $\REALS$, such that the inverse image
of any Borel set $B \in \BOREL$ belongs to $\SAF$, i. e. $X^{-1}(B) \in \SAF$, $\forall B \in \BOREL$.

Let's introduce the measure $\MUX$ for the measurable space $(\REALS,\BOREL)$ as:
for every subset $ B \in \BOREL$ define $\MUX(B) = \PM(X^{-1}(B))$. The measure $\MUX$ is called the measure {\em induced} by the random variable $X$, see theorem 3.1.3, \cite{CHUNG}.
Alternatively $\PM(X^{-1}(B))$ is also denoted as $\PM(X \in B)$, so that $\MUX(B) = \PM(X^{-1}(B)) = \PM(X \in B)$. 
Hence $(\REALS,\BOREL,\MUX)$ is a measure space, and as a matter of fact it is a probability space since 
$\MUX(\REALS) = \PM(X^{-1}(\REALS)) = \PM(X \in \REALS) = \PM(\SPO) = 1$. 

Define the distribution function of the random variable $X$ as 
$\DFX(x) = \PM( X \in (-\infty, x]) = \PM( X \le x)$. Note that some authors, e. g. \cite{LOEVE}, define the distribution function as the probability
of the inverse image of the open set $(\infty,x)$, i. e. $\DFX(x) = \PM( X < x)$,  instead of the set $(-\infty,x]$ as done here, following \cite{CHUNG}. 

The distribution function obtained
is a non-decreasing function, with:
\begin{equation}
\begin{aligned}
\nonumber
\lim_{ x \rightarrow -\infty} \DFX(x) = 0, 
\\
\lim_{ x \rightarrow +\infty} \DFX(x) = 1,
\end{aligned}
\end{equation}
which is continuous on the right at each $x \in \REALS$, and has a limit on the left at each $x \in \REALS$. As for the last property, if the distribution function is defined as in \cite{LOEVE}, the function is continuous on the left at each $x \in \REALS$, and has a limit on the right at each $x \in \REALS$.
Again, the definitions in \cite{CHUNG} are used here. Conversely, any function satisfying the preceding properties, is the distribution function of some random variable.

An additional most important property for what follows, is that the set of discontinuities of the distribution function, being an increasing function, is countable. 

A distribution function that is continuous everywhere is called a {\em continuous} distribution function.

Define $\STEP_{t}(x)$ as the step function $u:\REALS \rightarrow \REALS$, for which 
$\STEP_{t}(x) = 0$ if $x < t$, and $\STEP_{t}(x) = 1$ if $x \ge t$. 
Note that $\STEP_{t}(x) = \STEP_{0}(x-t)$.
The function $\STEP_{t}(x)$ is called the {\em point mass distribution at t}, \cite{CHUNG}. Note that $\STEP_{t}(x)$ is always continuous on the right.

A function that can be represented in the form 
$\sum_{l} \pi_{l} \; \STEP_{\xi_{l}}(x)$ where ${\xi_{l}}$ is a countable set of real numbers,
with $\xi_{l} \neq \xi_{n}$ if $ l \neq n$,
and where $\pi_{l} \ge 0$ for every $l$ and $\sum_{l} \pi_{l} = 1$, is called a {\em discrete} distribution function. Note that this function is well defined as a distribution function since it satisfies the previously mentioned properties. The corresponding random variable takes the value $\xi_{l}$ with probability $\pi_{l}$, and
for a value $x$ such that $x \neq \xi_{l}$, $\forall l$, it takes such a value with probability $0$.

A first result, see theorem 1.2.3, \cite{CHUNG}, is that every distribution function can be written as the convex combination of a discrete 
and a continuous distribution functions. Such a decomposition is unique.
That is, the distribution function $\DFX(x)$ may be decomposed as: 
$\DFX(x) = \rho \; \DFX^{d}(x) + (1 - \rho) \; \DFX^{c}(x)$, with $ 0 \le \rho \le 1$, where  $\DFX^{d}(x)$ and $\DFX^{c}(x)$ are respectively, the corresponding discrete and continuous distribution functions. As for $\DFX^{c}(x)$, since it is a bounded non-increasing function, then, it has non-negative derivative almost everywhere, but not necessarily this function is the integral of its derivative.

Define, \cite{CHUNG}, a function 
$F: \REALS \rightarrow \REALS$ as {\em absolutely continuous} in $(-\infty, +\infty)$ with respect to the {\em Lebesgue measure} $\lambda$, (i. e. $\lambda( (x, x^{\prime})) = x^{\prime} - x$, so that $\lambda(dt) = dt$), 
iff there exists a function $f$ in $L^{1}(\REALS, \BOREL, \lambda)$, such that for every $x < x^{\prime}$:
$F(x^{\prime}) - F(x) = \int_{x}^{x^{\prime}} f(t) \; dt$. Hence, an absolutely continuous distribution function $\DFX^{a}(x)$ is the integral of its derivative, which exists almost everywhere with respect to the Lebesgue measure $\lambda$, i. e. $\lambda$-ae. That derivative is defined as the {\em probability density function} of the random variable $X$ and will be denoted as $\PDFX(x)$, so that
$\DFX^{a}(x) = \int_{-\infty}^{x} \PDFX(t) \; dt$, and $\PDFX(x) = d \DFX^{a}(x) / dx$, $\lambda$-ae.

Define, \cite{CHUNG}, a distribution function 
$F: \REALS \rightarrow \REALS$ as {\em singular} iff it is not identically zero and $F^{\prime}$ exists and equals zero 
$\lambda$-ae, i. e. almost everywhere with respect to Lebesgue measure.

The main result is theorem 1.3.2, \cite{CHUNG}: every distribution function can be written as the convex combination of a discrete, 
a singular continuous, and an absolutely continuous distribution functions. This decomposition is unique.
That is $\DFX(x) = \rho_{1} \; \DFX^{d}(x) + \rho_{2} \; \DFX^{a}(x) + \rho_{3} \; \DFX^{s}(x)$, with $ 0 \le \rho_{i} \le 1$, for $i = 1, 2, 3$, and $\sum_{i=1}^{3} \rho_{i} = 1$, where  $\DFX^{d}(x)$, $\DFX^{a}(x)$, and $\DFX^{s}(x)$ are respectively, discrete, continuous, and singular continuous distribution functions.

Another basic important result, is given by theorem 3.2.2 , 
\cite{CHUNG}, see also \cite{LOEVE}, \cite{LIBROROJO}, \cite{H-S}. That is, let the random variable $X$ on $(\SPO,\SAF,\PM)$ induce the probability space 
$(\REALS, \BOREL, \MUX)$, and let $f: \REALS \rightarrow \REALS$ be Borel measurable, then:
\begin{equation}
\nonumber
\int_{\SPO} f(X(\omega)) \; \PM(d\omega) = \int_{\REALS} f(x) \; \MUX(dx) = \int_{-\infty}^{+\infty} f(x) \; d\DFX(x),
\end{equation}
provided that any of these integrals exists. Note that these integrals are different in nature: the first is a Lebesgue integral defined over the space $(\SPO,\SAF,\PM)$,
the second is a Lebesgue integral defined over the space $(\REALS, \BOREL, \MUX)$, while the third is a Lebesgue-Stieljes integral in the real line.

Let $s$ be a statement taking the values $TRUE$ or $FALSE$. Define the function taking real values $\UNOG(s)$ as $\UNOG(s) = 1$ if $s = TRUE$, while $\UNOG(s) = 0$ if $s = FALSE$.
For a given set $B \subseteq \REALS$, define $\UNOG_{B}(x) = \UNOG( x \in B)$, i. e. $\UNOG_{B}(x)$ is the characteristic or indicator function of the set $B$.
Abusing notation put $\UNOG_{\xi}(x)$ for $\UNOG_{\{\xi\}}(x)$ if $\xi$ is an isolated point. 
Define $\UNOG^{*}_{B}(x) = 1 - \UNOG_{B}(x)$ and put $\UNOG^{*}_{\xi}(x)$ for $\UNOG^{*}_{\{\xi\}}(x)$.

Let $B \in \BOREL$, and take $f(x) = \UNOG_{B}(x)$ in the previous result, then:
\begin{equation}
\nonumber
\int_{\REALS} \UNOG_{B}(x) \MUX(dx) = \int_{B} \MUX(dx) = \MUX(B) = \int_{-\infty}^{+\infty} \UNOG_{B} d\DFX(x).
\end{equation} 
This, together with
the definition of the distribution function, shows that there is a one to one correspondence between the induced measure and the distribution function for a given random variable.

Call a {\em discrete random variable}, a random variable whose distribution function is discrete, and call {\em absolutely continuous random variable}, a random variable whose distribution function is absolutely continuous. 

For the cases considered here, let's assume that the random variable $X$ has no singular component, that is, the distribution function 
of the random variable $X$ is a convex combination of a discrete and an absolutely continuous (with respect to Lebesgue measure) distribution functions. Call such a random variable a {\em real mixed states random variable}. 
Hence, the distribution function $\DFX^{m}$ of a real mixed states random variable,
has the form: $\DFX^{m}(x) = \rho \; \DFX^{d}(x) + (1 - \rho) \; \DFX^{a}(x) 
= \rho \; \sum_{l} \pi_l \; \STEP_{\xi_l}(x) + (1 - \rho) \; \int_{-\infty}^{x} \PDFX(t) \; dt$, with $ 0 \le \rho \le 1$. 
It is readily checked that $\MMUX$, the induced measure by the real mixed states random variable $X$, has the form:
$\MMUX(B) = \rho \; \sum_{l} \pi_l \; \UNOG_{B}(\xi_l) + (1 - \rho) \; \int_{B} \PDFX(t) \; dt$,
for each $B \in \BOREL$. 

In \cite{L-R} there is a study on the decomposition of the distribution function of real random variables in $\REALS^{n}$.
In particular it is shown that discontinuities in the distribution function occur in planes parallel to the axes.
The nature of discontinuities in $\REALS^{n}$, is far more complicated due to the possibility of probability mass concentration in
hyper-volumes of dimension less than $n$. It may happen that there is mass concentration in curves, surfaces, volumes, or hyper-volumes of dimension less than $n$, which are not presently discussed here. 

\subsection{Probability density function}\label{sec:PDF-REAL}

Here, the probability density for a real mixed states random variable will be analyzed.
As previously discussed, when a random variable $X$ has an absolutely continuous distribution function, then such a function is the integral of its derivative, which exists $\lambda$-ae. In such a case the derivative is called the {\em probability density function}. For a real mixed states random variable, its distribution function is a convex combination of a discrete and an absolutely continuous distribution functions with respect to the Lebesgue measure:
$\DFX^{m}(x) = \rho \; \DFX^{d}(x) + (1 - \rho) \; \DFX^{c}(x) 
= \rho \; \sum_{l} \pi_l \; \STEP_{\xi_l}(x) + (1 - \rho) \; \int_{-\infty}^{x} \PDFX(t) \; dt$, with $ 0 \le \rho \le 1 $. 
If $ \rho > 0 $, then the distribution function $\DFX^{m}(x)$ is not absolutely continuous with respect to the Lebesgue measure and then a probability density can not be defined. 
In this case two possible directions may be taken.

A first direction would be to consider a generalized probability density function $p^{m}_{X}(x)$ in the distributional sense, \cite{LAURENT-SCHWARTZ},
\cite{GELFAND},  loosely in the form:
$p^{m}_{X}(x) = \rho \; \sum_{l} \pi_{l} \; \delta(x - \xi_{l}) + (1 - \rho) \; p^{a}_{X}(x)$, $0 \le \rho \le 1$, where $\delta(x)$ is the Dirac delta "function" centered at $x = 0$, $\pi_{l}$ is the probability mass at $x = \xi_l$,
and $p^{a}(x)$ is the probability density of an absolutely continuous distribution function with respect to Lebesgue measure. 
Loosely one has $\UNOG_{B}(0) = \int_{B} \delta(x) dx$, and
$\STEP_{a}(x) = \STEP_{0}(x - a) = \int_{-\infty}^{x} \delta(t - a) \; dt$.
As is well known the Dirac "function" $\delta(x)$ is not a function in the ordinary sense, so that necessarily such an approach should be formally based on the 
theory of distributions, in either the view of \cite{LAURENT-SCHWARTZ}, or the view of \cite{GELFAND}.
This approach may also be loosely viewed as a mixture of a discrete and a continuous random variable. As a matter of fact, this intuitive approach was taken in \cite{SPIE-06}, considering the previous distributional density as the limit of a mixture of Gaussian random variables, one centered at zero with a a fixed very low variance $\sigma_0$, i. e. very small with respect to the variance of the second Gaussian component, so that the first Gaussian component could loosely be interpreted as a Dirac delta function. In \cite{ICIP-06}, \cite{SPIE-06} Markov Random Fields with mixed states random variables were used, where there appears a number of technical problems when dealing with the Dirac delta generalized function. In \cite{SPIE-06} to circumvent those difficulties and give a rather intuitive presentation of the theory, a mixture of Gaussian random variables was used to give a feeling of the results.
Formalizing the distributional approach seems too complicated when dealing with mixed states Markov Random Fields, so that to formalize these 
models a measure theoretic approach is proposed.

\subsection{Measure theoretic approach}\label{sec:measure}

As previously discussed, the second approach is measure theoretic. 
A formal approach is presented next, and theoretical results are given regarding the structure of real mixed states random variables, as to obtain general results that permit to construct models using these variables, such  as those
previously discussed in \cite{Bouthemy-1}, \cite{SPIE-06}, \cite{ICIP-06}, \cite{Hardouin-1}, \cite{Hardouin-2}, 
\cite{Gwenaelle-5}, as well as other estimation, classification, 
segmentation, detection, and filtering problems in more complicated situations presently under study. 

Recall that a measure $\nu$ is {\em absolutely continuous} with respect to a measure $\zeta$ both for the same measurable space, if for all measurable subsets $B$ such that $\zeta(B) = 0$ it is $\nu(B) = 0$.

Call a measure of the form $\nu(B) = \sum_{i} c_i \; \UNOG_{B}(d_i)$ with $c_i \ge 0$ for the measurable space $(\REALS, \BOREL)$, a {\em discrete measure}, whether $\sum_{i} c_i$, is finite or not, and a measure of the form
$\nu(B) = \int_{B} g \; d\lambda$, with $g: \REALS \rightarrow \REALS$ a fixed non-negative measurable function, $g \ge 0$, $\lambda$-ae, an {\em absolutely continuous measure with respect to the Lebesgue measure $\lambda$}, whether the integral is finite or not.
Let $f:\REALS \rightarrow \REALS$ be an arbitrary non-negative measurable function, then, if $\nu$ is a discrete measure:
$\int f d\nu = \sum_{i} c_{i} f(d_{i})$, while if $\nu$ is absolutely continuous w.r.t. $\lambda$ it is:
$\int f \; d\nu = \int f \; g\; d\lambda$. Also, if $\nu$ is of the form $\nu = \nu^{d} + \nu^{a}$, where $\nu^{d}$ and $\nu^{a}$ are respectively discrete and
absolutely continuous w.r.t. $\lambda$, then:
$\int f \; d\nu = \sum_{i} c_{i} \; f(d_{i}) + \int f \; g \; d\lambda$. 
It follows that this result is also true for $f \in L^{1}(\REALS, \BOREL, \nu)$.

Recall that a real mixed states random variable, is defined as a random variable with distribution function 
$\DFX^{m}(x) = \rho \; \DFX^{d}(x) + (1 - \rho) \; \DFX^{a}(x)$.

As previously discussed the induced measure generated by a real mixed states random variable $X$ takes the form:
\begin{equation}\label{RMMUX}
\MMUX(B) = \rho \; \MUX^{d}(B) + (1 - \rho) \; \MUX^{a}(B) = \rho \; \sum_{l} \pi_l \; \UNOG_{B}(\xi_l) + (1 - \rho) \; \int_{B} \PDFX(t) \; dt,
\end{equation}
for each $B \in \BOREL$, 
where both $\MUX^{d}$ and $\MUX^{a}$ are probability measures.
Hence, the Lebesgue integral of a non-negative measurable function $f: \REALS \rightarrow \REALS$ with respect to the measure space $(\REALS, \BOREL, \MUX)$ is:
$\int f(t) \; \MMUX(dt) = \rho \; \sum_{l} \pi_{l} \; f(\xi_{l}) + (1 - \rho) \; \int_{\REALS} f(t) \; \PDFX(t) \; dt$. It follows that this result is also true
for $f \in L^{1}(\REALS, \BOREL, \MMUX)$.

Let's introduce the measure $m(B)$, for each $B \in \BOREL$, for the measurable space $(\REALS,\BOREL)$ as: 
\begin{equation}\label{MB}
m(B) = \sum_{l} \UNOG_{B}(\xi_l) + \int_{B} dx = \MDX (B) + \lambda(B),
\end{equation} 
for each $B \in \BOREL$, where $\lambda(B)$ is the Lebesgue measure of $B$, and
$\MDX$ is the discrete measure $\MDX(B) = \sum_{l} \UNOG_{B}(\xi_l)$.

Next, let's show that the induced measure $\MMUX$ is absolutely continuous with respect to the measure $m$, that is $\MUX << m$.
The result to be proved is that if for $B \in \BOREL$ it is $m(B) = 0$, then $\MMUX(B) = 0$.
Suppose $B \in \BOREL$ is such that $m(B) = 0$, then $\UNOG_{B}(\xi_l) = 0$, for each $l$, so that $\xi_l \notin B$, for each $l$,
and then $\sum_{l} \pi_{l} \UNOG_{B}(\xi_l) = 0$.
Also, since $m(B) = 0$, then $\lambda(B) = 0$, but then, $\int_{B} \PDFX(t) \; dt =0$. Combining these two results $\MMUX(B) = 0$ is obtained.

Given the probability space $(\SPO, \SAF, \PM)$, the mixed states random variable $X$ induces the probability space $(\REALS, \BOREL, \MUX)$.
The main idea is to consider the measure space $(\REALS, \BOREL, m)$ instead of the standard Borel measure space $(\REALS, \BOREL, \lambda)$.
The main point is that when the probability space $(\REALS, \BOREL, \MUX)$, is referred with respect to the measure space $(\REALS, \BOREL, m)$, 
using the Radon-Nikodym theorem, \cite{H-S}, a generalized probability density function may be defined, that permits to handle simultaneously discrete and continuous valued random  variables.

The version of the Radon-Nikodym theorem used here is, \cite{H-S}: let $\nu$ be a finite measure on the measurable space 
$(\REALS, \BOREL)$, and $\zeta$ be a 
$\sigma$-finite measure for the same measurable space $(\REALS, \BOREL)$, then, if $\nu << \zeta$, i. e. $\nu$ is absolutely continuous with respect 
to $\zeta$, then, there exists a non-negative measurable $\zeta$-ae function $g: \REALS \rightarrow \REALS$, such that for all $B \in \BOREL$, it is
$\nu(B) = \int_{B} g \; d\zeta$. Note that the function $g$ is unique $\zeta$-ae.

Also, if $\nu(B) = \int_{B} g \; d\zeta$ is true for all $B \in \BOREL$, then, \cite{H-S}, if $f: \REALS \rightarrow \REALS$ is any non-negative $\nu$-ae measurable function, or  $f \in L^{1}(\REALS, \BOREL, \nu)$, then:
$\int f \; d\nu = \int f \;  g \; d\zeta$.

Call $A$ the countable set of all the $\xi_l$'s, i. e. $A = \{ \xi_l : \; l \in \mathbb{N}\}$, 
and call $\UNOG_{\xi_l}(x)$ the {\em point mass probability at $\xi_l$}.

Next, apply the Radon-Nikodym theorem to the induced measure $\MUX$, which is absolutely continuous 
with respect to the previously introduced measure $m$. Let's check that
the Radon-Nikodym derivative is:
\begin{equation}\label{RND}
\frac{d \MMUX}{d m} \equiv\MPDFX(x) = \rho \;  \sum_{l} \pi_{l} \; \UNOG_{\xi_l}(x) + (1 - \rho) \;  \UNOG^{*}_{A}(x) \; \PDFX(x). 
\end{equation}

To prove this result let's calculate, see equation \eqref{MB}:
\begin{equation}
\label{RNP}
\int_{B} \MPDFX \; dm= \int_{B} \MPDFX \; d \MDX + \int_{B} \MPDFX \; d \lambda.
\end{equation}
For the first integral:
\begin{equation}
\nonumber
\int_{B} \MPDFX \; d \MDX = \int_{B} (\rho \; \sum_{l} \pi_{l} \; \UNOG_{\xi_l}(x) + (1 - \rho) \; \UNOG^{*}_{A}(x) \; \PDFX(x)) \; \MDX(dx).
\end{equation}
For the first term, $\int_{B} \sum_{l} \pi_{l} \; \UNOG_{\xi_l}(x) \; \MDX(d x) = 
\int \UNOG_{B}(x) \; \sum_{l} \pi_{l} \; \UNOG_{\xi_l}(x) \; \MDX(dx) =  \sum_{j} \; \UNOG_{B}(	\xi_{j}) \; \sum_{l} \pi_{l} \; \UNOG_{\xi_l}(	\xi_{j}) =
 \sum_{j} \; \pi_{j} \; \UNOG_{B}(\xi_{j})$, since $ \UNOG_{\xi_l}(\xi_{j}) = 1$ if $ \xi_l = 	\xi_{j}$, and $\UNOG_{\xi_l}(\xi_{j}) = 0$ if $\xi_l \neq \xi_{j}$.
As for the second term:
\begin{equation}
\begin{aligned}
\nonumber
\int_{B} \; \UNOG^{*}_{A}(x) \; \PDFX(x) \; \MDX(dx) 
&= 
\int \; \UNOG_{B}(x) \; \UNOG^{*}_{A}(x) \; \PDFX(x) \; \MDX(dx)
\\
&= \sum_{l} \UNOG_{B}(	\xi_{l}) \; \UNOG^{*}_{A}(	\xi_{l}) \; \PDFX(	\xi_{l}) = 0,
\end{aligned}
\end{equation} 
since $\UNOG^{*}_{A}(	\xi_{l}) = 0$, for all $ l \in \mathbb{N}$.
Hence, $\int_{B} \MPDFX \; d \MDX = \rho \; \sum_{l} \pi_{l} \; \UNOG_{B}(	\xi_{l})$.

For the second integral in equation \eqref{RNP}:
\begin{equation}
\begin{aligned}
\nonumber
\int_{B} \MPDFX \; d\lambda 
&= \int_{B} (\rho \; \sum_{l} \pi_{l} \; \UNOG_{\xi_l}(x) + (1 - \rho) \; \UNOG^{*}_{A}(x) \; \PDFX(x)) \; d x
\\
&= \rho \; \sum_{l} \pi_{l} \; \int_{B} \; \UNOG_{\xi_l}(x) \; d x + (1 - \rho) \; \int_{B} \UNOG^{*}_{A}(x) \; \PDFX(x) \; d x, 
\end{aligned}
\end{equation}
applying Tonelli. 
But $ \lambda(\xi_l) = 0 $, for $l = 1, 2, \cdots$, since a single point has Lebesgue measure zero, hence the first term is zero. As for the second term, since
$\lambda(A) = 0$, because A is a countable collection of points, then:
\begin{equation} 
\nonumber
\int_{B} \UNOG^{*}_{A}(x) \; \PDFX(x) \; dx = \int_{B}  \PDFX(x) \; dx, 
\end{equation}
so that $\int_{B} \MPDFX \; d\lambda = (1 - \rho) \; \int_{B}  \PDFX(x) \; dx$.
Collecting all these results one obtains:
$\int_{B} \MPDFX \; dm = \rho \; \sum_{l} \pi_{l} \UNOG_{B}(	\xi_{l}) + (1 - \rho) \; \int_{B}  \PDFX(x) \; dx = \MMUX(B)$, 
for all $B \in \BOREL$, which is the desired result, see equation \eqref{RMMUX}.

\subsection{Label-real mixed states random variables}\label{sec:labels}

The label-real mixed states random variables are defined as to consider probability mass concentrated either in symbolic labels, 
as well as  values in $\REALS^{n}$ when these values belong to at most a countable subset of $\REALS^{n}$. 
If this is the case, all those values will belong to the set $\ELE$.
Let $\ELE = \{ \ell_1, \ell_2, \cdots \}$ be a countable set of symbolic labels and eventually at most a countable number of values in $\REALS^{n}$.
Define the mixed states space as $\EME = \ELE \cup \REALS^{n}$, with $n \ge 1$. Define $\SAM$ as the collections of subsets of $\EME$, such that
$\SAM = 2^{\ELE} \cup \BOREL(\REALS^{n})$, where $2^{\ELE}$, is the power set of $\ELE$, i. e. the collection of all the subsets of $\ELE$, and
$\BOREL(\REALS^{n})$ is the Borel $\sigma$-algebra for $\REALS^{n}$, i. e. the minimal $\sigma$-algebra generated by, say, the open sets of $\REALS^{n}$.
Redefine $\lambda$ as the Lebesgue measure for $(\REALS^{n}, \BOREL(\REALS^{n}))$, i. e. the measure that assigns to a hyper-cube in $\REALS^{n}$
its volume, given by the product of the length of the sides of the hyper-cube, so that $\lambda(d \X) = d \X$.

If $M \subseteq \EME$, define $\complement(M)$, the complement of $M$ as the set $\complement(M) = \EME \setminus M$, i. e. the elements in $\EME$ 
not in $M$. Note that any set $M \in \SAM$, can be decomposed as $M = D \cup B$, with 
$D \in 2^{\ELE}$ and $B \in \BOREL(\REALS^{n})$. Such decomposition is unique if $\varnothing = 2^{\ELE} \cap \BOREL(\REALS^{n})$, i. e. the
set $\ELE$ consists only of symbolic labels. 

Clearly $\SAM$ is a $\sigma$-algebra since: i) $\varnothing \in \SAM$, ii) if $M \in \SAM$, then there exists
$D \in 2^{\ELE}$ and $B \in \BOREL(\REALS^{n})$, with $\varnothing = D \cap B$ and $\varnothing = B \cap \ELE$, such that $M = D \cup B$, so that
$\complement(B) = (\ELE \setminus D) \cup (\REALS^{n} \setminus B)$, and then $\complement(M) \in \SAM$, iii) if $\{M_i\}_{i \in \mathbb{N}}$ is a countable sequence of sets $M_i \in \SAM$, $i = 1, 2, \cdots$, then there exist $D_i \in 2^{\ELE}$, and $B_i \in \BOREL(\REALS^{n})$, $i = 1, 2, \cdots$, such that
$M_i = D_i \cup B_i$, $i = 1, 2, \cdots$, but since $\cup_{i \in \mathbb{N}} D_i \in 2^{\ELE}$ and $\cup_{i \in \mathbb{N}} B_i \in \BOREL(\REALS^{n})$,
then $M \in \SAM$.

Hence, a measurable space $(\EME, \SAM)$ was constructed. Following \cite{CHUNG}, given the probability space $(\SPO, \SAF, \PM)$, define the function 
$\XRV: \SPO \rightarrow \EME$
as a mixed states random variable if $\XRV^{-1}(M) \in \SAF$, for all $M \in \SAM$, 
i. e. the inverse image under $\XRV$ of any set in the $\sigma$-algebra $\SAM$
is in the $\sigma$-algebra $\SAF$.

The main idea now, is to construct an induced measure $\MMUXV$ for the random variable $\XRV$ in the measurable space $(\EME, \SAM)$. 
Since now symbolic labels are present, which may not have
any algebraic structure, a distribution function can not be defined to characterize the random variable. Hence, a possibility is to proceed directly to define the
measure $\MMUXV(M)$, for each $M \in \SAM$, as, see equation \eqref{RMMUX}:
\begin{equation}\label{eq:1001}
\begin{aligned}
 \MMUXV(M) &= \rho \; \sum_{l} \pi_l \; \UNOG_{M}(\ell_l) + (1 - \rho) \; \int_{M \setminus \ELE} \PDFXV(\X) \; d \X
\\
 &= \rho \; \MUXV^{d}(M) + (1 - \rho) \; \MUXV^{a}(M), 
\end{aligned}
\end{equation}
where $\MUXV^{d}(M) = \sum_{l} \pi_l \; \UNOG_{M}(\ell_l)$ and $ \MUXV^{a}(M) = \int_{M \setminus \ELE} \PDFXV(\X) \; d \X$, 
with $ 0 \le \rho \le 1$, $0 \le \pi_l \le 1$, $ l = 1, 2, \cdots$, 
$\sum_{l} \pi_l = 1$, and
 $\PDFXV(\X)$, the probability density with respect to Lebesgue measure of some real multidimensional standard random variable.
Note that $\MMUXV(\EME) = 1$, so that $(\EME, \SAM, \MMUXV)$ is a probability space. Also, note that $\int_{M \setminus \ELE} \PDFXV(\X) \; d \X \equiv
\int_{M \setminus \ELE} \PDFXV \; d \lambda$, is a standard Lebesgue integral in $(\REALS^{n}, \BOREL(\REALS^{n}))$, since the integral is evaluated
in the set $M \setminus \ELE \subseteq \REALS^{n}$, with $M \setminus \ELE \in \BOREL(\REALS^{n})$.

Clearly if $\rho > 0$ then $\MMUXV$ is not absolutely continuous with respect to the Lebesgue measure.
Following the same procedure as before, redefine the reference measure $m$ for the measurable space $(\EME, \SAM)$ as, see equation \eqref{MB}:
\begin{equation}\label{LMB}
m(M) = \sum_{l} \UNOG_{M}(\ell_l) + \int_{M \setminus \ELE} d \X = \MDX (M) + \tilde{\lambda}(M),
\end{equation} 
for each $M \in \SAM$, where now $\MDX$ is the discrete measure $\MDX(M) = \sum_{l} \UNOG_{M}(\ell_l)$
and where $\tilde{\lambda}(M) = \int_{M \setminus \ELE} d \X$.  Note that $\tilde{\lambda}(M)$ is not the Lebesgue measure of $M$, since
$M \notin \BOREL(\REALS^{n})$ in general, but is close in the
sense that $\tilde{\lambda}(M) = \lambda(M)$ if $M \in \BOREL(\REALS^{n})$,
because, since $\tilde{\lambda}(M) = \lambda (M \setminus \ELE)$, and
$M \subseteq \REALS^{n}$ has at most a countable number of elements in $\ELE$ which have Lebesgue measure zero, 
then $\lambda(M \setminus \ELE) = \lambda(M)$, so that
$\tilde{\lambda}(M) = \lambda(M)$.

Also, $\int_{M} f \; d\tilde{\lambda} = \int_{M \setminus \ELE} f \; d \lambda$, for all non-negative measurable $f$, and then for all 
$f \in L^{1}(\EME, \SAM, \tilde{\lambda})$. To show this result, start first, as usual, with characteristic functions.
Hence, let $f = \UNOG_M$, with $M \in \SAM$. Then, $\int \UNOG_M \; d\tilde{\lambda} = \tilde{\lambda}(M) =
\int_{M \setminus \ELE} \; d \lambda = \int_{\REALS^{n} \setminus \ELE} \UNOG_M\; d \lambda$. Next, \cite{LIBROROJO}, \cite{HALMOS}, \cite{H-S}, proceed to simple functions, and then to non-negative functions using the monotone convergence theorem on both sides of the equality. Finally, the result
for $f$ non-negative, or for $f \in L^{1}(\EME, \SAM, \tilde{\lambda})$, is obtained as $\int f \; d\tilde{\lambda} = \int_{\REALS^{n} \setminus \ELE} f \; d \lambda$.
Hence for non-negative measurable $f$, or for $f \in L^{1}(\EME, \SAM, \tilde{\lambda})$, and arbitrary $M \in \EME$:
$\int \UNOG_M \; f \; d\tilde{\lambda} = \int_{\REALS^{n} \setminus \ELE} \UNOG_M\; f \; d \lambda$, so that 
$\int_{M} f \; d\tilde{\lambda} = \int_{M \setminus \ELE} f \; d \lambda$.

As before, let's proceed to show that $\MMUXV$ is absolutely continuous with respect to $m$, i. e. $\MMUXV << m$.
If $m(M) = 0$, then $\int_{M \setminus \ELE} d \X = 0$ and $\sum_{l} \UNOG_{M}(\ell_l) = 0$. Hence, $\UNOG_{M}(\ell_l) = 0$, for $l = 1, 2 , \cdots,$ so that
$\ell_l \notin M$, for $l = 1, 2 , \cdots$, hence $\varnothing = M \cap \ELE$, and then $M \setminus \ELE = M$, and $M \in \BOREL(\REALS^{n})$.
Then $\int_{M \setminus \ELE} d \X = \int_{M} d \X = \lambda(M) = 0$. Hence $M$ is Borel measurable, and has Lebesgue measure zero, so that
$\int_{M \setminus \ELE} \PDFXV({\X}) \; d{\X} = \int_{M} \PDFXV({\X}) \; d{\X} = 0$. Also, since $\ell_l \notin M$, for $l = 1, 2 , \cdots$, then  
$\sum_{l} \pi_l \; \UNOG_{M}(\ell_l) = 0$, so that $\MMUXV(M) = 0$. 

Hence, since $m$ is a $\sigma$-finite measure, and $\MMUXV$ is a finite measure, as a matter of fact a probability measure,
consider, as it was previously done, the Radon-Nikodym formalism, \cite{H-S}. As in the previous section, let's check that
the Radon-Nikodym derivative is, see equation \eqref{RND}:
\begin{equation}\label{RNDV}
\frac{d \MMUXV}{d m} \equiv \MPDFXV(\X) = \rho \;  \sum_{l} \pi_{l} \; \UNOG_{\ell_l}(\X) + (1 - \rho) \;  \UNOG^{*}_{\ELE}(\X) \; \PDFXV(w(\X)),
\end{equation}
where the function $w: \EME \rightarrow \REALS^{n}$ is such that $w(\X) = \X$ if $\X \notin \ELE$, while $w(\X) = \X_{\bf R}$ if $\X \in \ELE$,
where $\X_{\bf R} \in \REALS^{n}$ is a fixed value in the domain of the function $\PDFXV$. Note that which value is chosen for $\X_{\bf R}$ is of no
importance, as long as is a valid value for $\PDFXV$, because of the factor $\UNOG^{*}_{\ELE}(\X) = 0$ if $\X \in \ELE$, resulting that the value of the second term in equation \eqref{RNDV} is not altered by the choice $\X_{\bf R}$.

To prove the result given by equation \eqref{RNDV} let's calculate, see equation \eqref{LMB}:
\begin{equation}
\label{RNPV}
\int_{M} \MPDFXV \; dm= \int_{M} \MPDFXV \; d \MDX + \int_{M} \MPDFXV \; d \tilde{\lambda}.
\end{equation}
For the first integral:
\begin{equation}
\nonumber
\int_{M} \MPDFXV \; d \MDX = \int_{M} \left(\rho \; \sum_{l} \pi_{l} \; \UNOG_{\ell_l}(\X) + (1 - \rho) \; \UNOG^{*}_{\ELE}(\X) \; \PDFXV(w(\X))\right) \; \MDX(d \X).
\end{equation}
For the first term:
\begin{equation}
\begin{aligned}
\nonumber
\int_{M} \sum_{l} \pi_{l} \; \UNOG_{\ell_l}( \X) \; \MDX( d \X) 
&= \int \UNOG_{M}(\X) \; \sum_{l} \pi_{l} \; \UNOG_{\ell_l}(\X) \; \MDX(d \X)
\\
&= \sum_{j} \; \UNOG_{M}(\ell_{j}) \; \sum_{l} \pi_{l} \; \UNOG_{\ell_l}(\ell_{j})
= \sum_{j} \; \pi_{j} \; \UNOG_{M}(\ell_{j}), 
\end{aligned}
\end{equation}
since $ \UNOG_{\ell_l}(\ell_{j}) = 1$ if $ \ell_l = \ell_{j}$,
and $\UNOG_{\ell_l}(\ell_{j}) = 0$ if $\ell_l \neq \ell_{j}$.
As for the second term:
\begin{equation}
\begin{aligned}
\nonumber
\int_{M} \; \UNOG^{*}_{\ELE}(\X) \; \PDFXV(w(\X)) \; \MDX(d \X) 
&= 
\int \; \UNOG_{M}(\X) \; \UNOG^{*}_{\ELE}(\X) \; \PDFXV(w(\X)) \; \MDX(d \X)
\\
&= \sum_{l} \UNOG_{M}(	\ell_{l}) \; \UNOG^{*}_{\ELE}(\ell_{l}) \; \PDFXV(\X_{\bf R}) = 0,
\end{aligned}
\end{equation}
since $\UNOG^{*}_{\ELE}(\ell_{l}) = 0$, for all $l \in \ELE$. 
Hence, $\int_{M} \MPDFXV \; d \MDX = \rho \; \sum_{l} \pi_{l} \; \UNOG_{M}(\ell_{l})$.

For the second integral in equation \eqref{RNPV}:
\begin{equation}
\begin{aligned}
\nonumber
\int_{M} \MPDFXV \; d \tilde{\lambda} 
&= \int_{M \setminus \ELE} \MPDFXV \; d {\lambda} 
\\
&= \int_{M \setminus \ELE} (\rho \; \sum_{l} \pi_{l} \; \UNOG_{\ell_l}(\X) + (1 - \rho) \; \UNOG^{*}_{\ELE}(\X) \; \PDFXV(w(\X))) \; d \X
\end{aligned}
\end{equation}
\begin{equation}
\begin{aligned}
\nonumber
= \rho \; \sum_{l} \pi_{l} \; \int_{M \setminus \ELE} \; \UNOG_{\ell_l}(\X) \; d \X 
+ (1 - \rho) \; \int_{M \setminus \ELE} \UNOG^{*}_{\ELE}(\X) \; \PDFXV(w(\X)) \; d \X,
\end{aligned}
\end{equation}
applying Tonelli.
The first term is zero since for each $l = 1, 2, \cdots$, it is $\UNOG_{\ell_l}(\X) = 0$, for $\X \in M \setminus \ELE$.
As for the second term, since
$\UNOG^{*}_{\ELE}(x) = 1$, for $ \X \in M \setminus \ELE$, and since $w(\X) = \X$ for $\X \in M \setminus \ELE$, then:
\begin{equation}
\nonumber
\int_{M \setminus \ELE} \UNOG^{*}_{\ELE}(x) \; \PDFXV(w(\X)) \; d \X = \int_{M \setminus \ELE}  \PDFXV(\X) \; d \X,
\end{equation}
so that $\int_{M} \MPDFXV \; d \tilde{\lambda} = (1 - \rho) \; \int_{M \setminus \ELE}  \PDFXV(\X) \; d \X$.
Collecting all these results one obtains:
$\int_{M} \MPDFXV \; dm = \rho \; \sum_{l} \pi_{l} \UNOG_{M}(\ell_{l}) + (1 - \rho) \; \int_{M \setminus \ELE}  \PDFXV(\X) \; d \X = \MMUXV(B)$, 
for all $B \in \BOREL(\REALS^{n})$, which is the desired result, see equation \eqref{eq:1001}.

\section{Real Mixed States Markov Random Fields}\label{sec:R-MS-MRF}

\subsection{Some definitions}\label{sec:defs}

Let $(\REALS^{N}, \BOREL(\REALS^{n}))$ be the measurable space where the space is $\REALS^{N}$, and $\BOREL(\REALS^{n})$ is the $\sigma$-algebra generated say, by the open sets in $\REALS^{N}$. Let $A$ be a finite set. Define $\#(A)$ as the number of elements in $A$.
Consider now the following spaces $\mathcal{D}_{i}$, where for each $i$, 
each space $\mathcal{D}_{i}$ might be a discrete space, or $\mathbb{R}$, or $\mathbb{R}^{n_{i}}$
or $\mathbb{C}^{n_{i}}$, for some $n_i \ge 1$, 
or the union of a discrete space and $\mathbb{C}^{n_{i}}$ or $\mathbb{R}^{n_{i}}$, i. e.
a mixed states variable. 
Here, each $\mathcal{D}_{i}$ will be considered as the union
of a single discrete value, where probability mass will be concentrated, and $\mathbb{R}^{n_{i}}$, with 
eventually a different dimension $n_{i}$ for each $i$.
Let $N$ be a finite positive integer and define 
$\mathcal{D} = \mathcal{D}_{1} \times \mathcal{D}_{2} \times \cdots \times \mathcal{D}_{N}$.
Let $f$ denote a function from $\mathcal{D}$ to $\mathbb{R}$. The function $f$ is a function
with $N$ arguments (called also "sites") for which the $i$-th argument takes 
values $\XK \in \mathcal{D}_{i}$. 

Let $\R \in \mathcal{D}$ be a fixed vector such that
$f(\R) = 0$. Call that vector a reference vector, or "ground" vector.
Each $\RK \in \mathcal{D}_{i}$ will be called the "ground" value of the $i$-th argument.
Equivalently, let's say that the $i$-th component of a vector $\X \in \mathcal{D}$ is
"grounded" if $\XK = \RK$.

Let $S$ be the set $S = \{1, 2, \cdots, N\}$. For a given $A \subseteq S$ define $\XA$ as a vector built as the concatenation of 
the vectors $\XK$ for $i  \in A$. The superscript $^{c}$ when applied to a set, denotes the  set complementing operation. In this section set complementing is considered with respect to the set $S$.

Define $\XKC$ as a vector built as the concatenation of all the $\XJ$ for $j \neq i$, i.e.
the vector $\XKC$ does not contain $\XK$. Note that $\XKC = \X_{\{i\}^{c}}$.

For a given $A \subseteq S$ define $g_{A}(\X)$ as a vector of the same dimension
as $\X$ such that the $i$-th component is $[g_{A}(\X)]_{i} = \XK$ if $i \in A$ while $[g_{A}(\X)]_{i} = \RK$
if $i \notin A$. That is, the values corresponding to sub-indexes in $A$ are preserved while those corresponding to sub-indexes not
in $A$ are "grounded". Call the function $g_{A}(\X)$ the "grounding" function.

\subsection{The Hammersley-Clifford Theorem}
\label{sec:HC-theorem}

\begin{theorem}[\cite{BESAG-1}] \label{th:BHC-1}
Let $f$ be a function from $\mathcal{D}$ to any arbitrary group. In particular, that group may be taken as $\REALS^{n}$, for some $n \ge 1$.
Call $S = \{1,2, \cdots, N \}$. Given $\R \in \mathcal{D}$ such that $f(\R) = 0 $,
then there is a decomposition:
\begin{equation}\label{eq1}
f(\X) = \sum_{A \subseteq S} f_{A}(\XA)
\end{equation}
satisfying:
\\
\noindent
A) Each $f_{A}(\X_{A})$ is a function that only depends on the arguments whose sub-indexes belong
to the set A, ordered lexicographically, both, the sub-indexes in $A$ and the arguments in $f$, so that
there is a one to one correspondence between the elements of $A$ and the arguments of $f$.
That is, any change in the value of an argument not in $A$ does not alter the value
of $f_{A}(\X_{A})$.
\\
\noindent
B) If for at least some $i \in A$ it is $\XK =\RK$, then $f_{A}(\X_{A}) = 0$
\\
\noindent
C) For the given $\R$, the decomposition given in equation \eqref{eq1} is unique.

\end{theorem}

Note that as for the function $f$ from $\mathcal{D}$ to $\mathcal{G}$, it is only required that $\mathcal{G}$ has a group structure, i. e. the existence of a binary operation with group structure.

\begin{theorem} [Hammersley-Clifford theorem, \cite{BESAG-1}, \cite{K-S}] \label{th:HC-1}
~~
Let $ \; p(\X)$ be a generalized probability density with respect to some measure $m$,
such that $p(\X) \neq 0, \quad \forall \X \in \mathcal{D}$. 
Fix $\R \in \mathcal{D}$, then,
there exists a unique collection of functions $Q_{A}(\X_{A})$ 
which take the zero value if any of its arguments is grounded, such that:
\\
\noindent
\begin{equation}\label{eq2}
\LOG \frac{p(\X)}{p(\R)} = \CHC,
\end{equation}
and:
\begin{equation}\label{eq3}
\LOG \frac{p(\XK \vert \XKC)}{p(\RK \vert \XKC)} = \CHCC.
\end{equation}
\end{theorem}

\section{A Decomposition Theorem}\label{sec:theorem}

\subsection{Preliminaries}\label{subsec:theo-prelim}

In this section real mixed states Markov Random Fields will be considered. 
Following Chapter VI, \cite{H-S}, define the product measurable space $(\mathcal{D}, \mathcal{M})$, where $\mathcal{D}$ is the cartesian product $\mathcal{D} = \mathcal{D}_{1} \times \mathcal{D}_{2} \cdots \mathcal{D}_{N}$, and
$\mathcal{M}$ is the product $\sigma$-algebra generated by $\mathcal{M}_{1} \times \mathcal{M}_{2} \times \cdots \times \mathcal{M}_{N}$. Consider the probability space $(\mathcal{D}, \mathcal{M}, m)$,
where the measure $m$ is given by 
\begin{equation}\label{MEASURE}
m = \prod_{i = 1}^{N} m_{i},
\end{equation}
for rectangles, where each measure $m_i$ is of the
form $m_i = \UNOG_{\RK}+ \lambda_i$, where $\RK \in \REALS^{n_i}$, and where $\lambda_i$ is the Lebesgue measure
for the measurable space $(\REALS^{n_i}, \BOREL(\REALS^{n_i})$. Note that $m(\mathcal{D}) = 1$. Formalizing the product space is a rather
long issue. As said in \cite{H-S}, p. 379: "The subject exhibits several technicalities, which can be an annoyance or a source of fascination: depending upon one's point of view". Though this author belongs to the second group, the reader is referred to Chapter VI, \cite{H-S}, to raise the point on the construction of the product space. It can be checked
that there are no surprises and that equivalent results are obtained for the real mixed states random fields introduced here. In particular both
Tonelli's and Fubini's theorems are valid. Define as usual, the marginal densities and the conditional densities.

Define a {\em real mixed states random field},  as a collection of random variables $\bf X$ described by a joint generalized probability density function $\MRFX(\X)$ with respect to the measure $m = \prod_{i = 1}^{N} m_{i}$.  Note that the the function $\MRFX(\X)$ is a function from $\mathcal{D} \rightarrow \REALS$ that integrates to $1$ with respect to the measure $m$. Following \cite{K-S}, define a {\em real mixed states Markov random field}, denoted as ms-MRF, a real mixed states random field such that $\MRFX(\X) \neq 0$, for all $\X \in \mathcal{D}$.

In what follows drop the sub-index denoting the random variable $\bf X$  in $\MRFX(\X)$ and denote it as $\PMSX$, whenever there is no confusion in notation.

Let $\PKMS$ be a collection, for $i=1,2,\cdots,N$, of conditional probability density functions corresponding to random variables taking mixed values, i.e.
pdf's with respect to the measures $m_i(d \XK) = d (\XK + u_{\RK}(\XK))$, where $u_{\RK}(\XK)$ is a step function at $\XK = \RK$, as previously defined. 
Random variables whose pdf's are taken with respect to these measures $m_i$, will
be called real mixed states random variables and the corresponding densities will be called ms-pdf's
(mixed states probability density functions). Whenever a ms-pdf is integrated, it will
be done using a Lebesgue-Stieljes integral with respect to the corresponding measure $d (\XK + u_{\RK}(\XK))$ as previously discussed.

Recall:
\begin{equation}\label{def1}
\quad \DK \; = \; \left\{
\begin{aligned}
&1 \quad {\rm if} \quad \XK \; = \; \RK, \quad
\\
&0 \quad {\rm if} \quad \XK \; \neq \; \RK, \quad
\end{aligned}
\right.
\end{equation}
\noindent
and:
\begin{equation}\label{def2}
\DKS = 1 -\DK.
\end{equation}
Define:
\begin{equation}\label{def3p}
\QK = \mathcal{P}( \XK = \RK \; | \; \XKC),
\end{equation}
and:
\begin{equation}\label{def3}
\QKS = 1 -\QK.
\end{equation}

In this section let's determine the joint probability 
density function for which the conditional ms-pdf's with respect to the measures $m_i$ take the form:
\begin{equation}
\label{th1}
\begin{aligned}
\PKMS &= \QK \; \DK \; + \; \QKS \; \DKS \; \PCK,
\\
& \qquad \qquad \forall i \; = \; 1,2 \cdots, N.
\end{aligned}
\end{equation}
\noindent
These conditional pdf's can not be taken arbitrarily,
except when they are all independent, because they have to satisfy the Hammersley-Clifford theorem
in the form of theorem \ref{th:HC-1} in the previous section.
When the $\XK$'s are all independent, the conditional ms-pdf's 
may be taken arbitrarily. 
For this case the ms-pdf is simply the product of all the conditional ms-pdf's.
Additionally, the $\QK$'s with $0 \leq \QK \leq 1$, and the functions
$\PCK = \PC(\XK)$, may all be taken functionally independent from one another, in that case.

When the conditional ms-pdf's effectively depend on the neighbors, there appear additional
constraints that these functions must satisfy because of the Hammersley-Clifford theorem.

\subsection{The theorem}\label{thetheorem}

Let's say that $\PCX$ is a MRF {\it by itself} if it has the form given by equations 
\eqref{eq2}, \eqref{eq3} satisfying theorem \ref{th:HC-1},
so that the Markov Random Field is given by $\PCX$, which corresponds to the pdf of some continuous distribution function with respect
to Lebesgue measure. That is, if $\PCX$ is considered alone, as a pdf with respect to Lebesgue measure, then $\PCX$ is a MRF.

\begin{theorem}\label{th:THE-THEOREM}
If $p^{a}({\X})$ is by itself a Markov random field with respect to the 
Lebesgue measure for $(\REALS^{N}, \BOREL(\REALS^{n}))$, then, the potentials of the Gibbs representation, \cite{K-S}, of the joint ms-pdf, 
decompose in a discrete related part and a continuous related part,
from which all the conditions and constraints required  by the Hammersley-Clifford theorem are obtained.
\end{theorem}

Proof. From equation \eqref{th1}:
\begin{equation}
\label{th2}
\begin{aligned}
\LPKMS &= 
\DK \LQK + \; \DKS \LQKS
\\
& \quad \quad  + \; \DKS \LPCK
\\
&= \DK \LQK + \DKS \LQKS 
\\
& \quad \quad + \DKS \LFCK + \DKS \LPCRK
\\
&= \DK \LQK + \DKS \LQKS 
\\
& \quad \quad + \LFCK + \DKS \LPCRK
\end{aligned}
\end{equation}
\noindent
The last equality is true since $ \LFCK $ is zero for $ \XK = \RK $, and $\DKS = 1$ 
for $\XK \neq \RK$.
\\
\noindent
Also, since $\PCX$ is a MRF by itself:
\begin{equation}
\label{th3}
\LFPC \equiv \VCX \equiv \FCX = \CC,
\end{equation}
where to conform with the usual terminology, \cite{K-S}, $\VCX$ is called the energy of the MRF, and 
$\{Q^{a}_{A}(\X_{A})\}_{A \subseteq S}$ is the collection of the potentials of the MRF,
and:
\begin{equation}\label{th4}
\LFCK = \CCK.
\end{equation}
\noindent
Hence:
\begin{equation}\label{th5}
\begin{aligned}
&\LPKMS = \DK \LQK \; +
\\
&+ \DKS \LQSP
 + \CCK.
\end{aligned}
\end{equation}
\noindent
For $\XK = \RK $ use \eqref{th2} and \eqref{th3} to obtain:
\begin{equation}
\label{th6}
\LPRKMS = \LQK,
\end{equation}
so that:
\begin{equation}\label{th7}
\LFPKMS = \DKS  \LFQP  + \LFCK.
\end{equation}
\\
\noindent
Apply theorem \ref{th:HC-1} to $\PMSX$
to obtain:
\begin{equation}\label{th8}
\LFPMS \equiv \VMX \equiv \FMSX = \CMS,
\end{equation}
where, as before, $\VMX$ is the energy of the mixed states MRF, and the collection $\{Q^{m}_{A}(\X_{A})\}_{A \subseteq S}$ is the collection of potentials
of the mixed states MRF,
and:
\begin{equation}\label{th9}
\LFPKMS = \CMSK.
\end{equation}
\noindent
Applying \eqref{th4} and \eqref{th9} in \eqref{th7} it results:
\begin{equation}\label{th10}
\CMSK = \DKS  \LFQP + \CCK.
\end{equation}
Call:
\begin{equation}\label{th11}
\AK = \AKDEF,
\end{equation}
and define
\begin{equation}\label{th12}
\begin{aligned}
\GKXC = \LFQP -\AK,
\end{aligned}
\end{equation}
so that:
\begin{equation}\label{th13}
\begin{aligned}
\GKR = 0.
\end{aligned}
\end{equation}
\\
\noindent
Observe that for each $i$, each function $\GKXC$ may eventually depend 
on all the variables except the variable $\XK$.
Because of the $N$ equations \eqref{th13}, apply theorem \ref{th:BHC-1} to each of the $N$
functions in \eqref{th12} for $i=1,2,\cdots,N$, to obtain:
\begin{equation}\label{th14}
\begin{aligned}
\GKXC = \CGK, \quad {\rm for} \quad i = 1,2,\cdots,N.
\end{aligned}
\end{equation}
\noindent
Hence:
\begin{equation}\label{th15}
\begin{aligned}
\CMSK &= \DKS  \AK + \DKS \CGK 
\\
&+ \CCK, \quad {\rm for} \quad i = 1,2,\cdots,N.
\end{aligned}
\end{equation}

\subsubsection{First order cliques}

The first order cliques correspond to:
\begin{equation}\label{th201}
A \; = \; \{i\},  \quad {\rm for} \quad i = 1,2,\cdots,N.
\end{equation}
Instead of $\X$ use the grounding function $g_{\{i\}}(\X)$, so that $\XKC = \RKC$, in each of 
the $N$ equations given in \eqref{th15}, to obtain:
\begin{equation}\label{th203}
\begin{aligned}
\FOCKMS =  \AK \DKS + \FCCK,
\end{aligned}
\end{equation}
Hence, all the first order cliques for the mixed states distribution were obtained.

\subsubsection{Second order cliques}

The second order cliques correspond to:
\begin{equation}\label{th301}
A \; = \; \{i,j\},  \quad {\rm for} \quad 1 \leq i < j \leq N.
\end{equation}
Subtract from \eqref{th15} all the first order cliques obtained in \eqref{th203}. Next,
instead of $\X$ use the grounding function $g_{\{i,j\}}(\X)$, with $i < j$, to obtain:
\begin{equation}\label{th302}
\begin{aligned}
\FMS_{\{i,j\}}(\XK,\XJ) &= \DKS \BGK_{\{j\}}(\XJ) + \FC_{\{i,j\}}(\XK,\XJ)
\\
&= \DJS \BGJ_{\{i\}}(\XK) + \FC_{\{j,i\}}(\XJ,\XK).
\end{aligned}
\end{equation}
Hence:
\begin{equation}\label{th303}
\begin{aligned}
\DKS \BGK_{\{j\}}(\XJ) = \DJS \BGJ_{\{i\}}(\XK), \qquad \forall j \neq i, \quad 1 \leq i < j \leq N. 
\end{aligned}
\end{equation}
Fix a value $\XK^{\dagger}$ for $\XK$ such that $\XK^{\dagger} \neq \RK$, 
then $\DELTA_{\RK}^{*}( \XK^{\dagger}) = 1$. 
Call $\beta_{i,j} \; = \; \BGJ_{\{i\}}(\XK^{\dagger})$,
then:
\begin{equation}\label{th303-1}  
\BGK_{\{j\}}(\XJ) = \beta_{i,j}  \DJS.
\end{equation} 
Proceed analogously with $\XJ$ to obtain:
$\BGJ_{\{i\}}(\XK) = \beta_{j,i}  \DKS $. Hence, from \eqref{th303}:
\begin{equation}\label{th304}
\begin{aligned}
\DKS \DJS \beta_{i,j} = \DJS \DKS \beta_{j,i},
\end{aligned}
\end{equation}
so that:
\begin{equation}\label{th305}
\begin{aligned}
\beta_{i,j} =  \beta_{j,i}.
\end{aligned}
\end{equation}
Then, from \eqref{th302} the second order cliques may be obtained as:
\begin{equation}\label{th306}
\begin{aligned}
\FMS_{\{i,j\}}(\XK,\XJ) &= \; \beta_{i,j} \; \DKS \DJS + \FC_{\{i,j\}}(\XK,\XJ)
\end{aligned}
\end{equation}
\noindent
Hence, all the second order cliques for the mixed states distribution were obtained.

\subsubsection{Third and higher order cliques}

Proceed analogously for the higher order cliques. 
For clarity consider the third order cliques.
The third order cliques correspond to sets of the form:
\begin{equation}\label{th401}
A \; = \; \{i,j,k\},  \quad {\rm for} \quad 1 \leq i < j < k \leq N.
\end{equation}
Substract all the first  and second order cliques from \eqref{th15} using \eqref{th203}
and \eqref{th306}. Next,
instead of $\X$ use the grounding function $g_{\{i,j,k\}}(\X)$, with $i < j < k$, in the $i$-th, $j$-th,
and $k$-th equations to obtain:
\begin{equation}\label{th402}
\begin{aligned}
\FMS_{\{i,j,k\}}(\XK,\XJ, \XM) &= \DKS \BGK_{\{j,k\}}(\XJ, \XM) + \FC_{\{i,j,k\}}(\XK,\XJ,\XM)
\\
&= \DJS \BGJ_{\{i,k\}}(\XK,\XM) + \FC_{\{j,i,k\}}(\XJ,\XK,\XM)
\\
&= \DMS \BGM_{\{i,j\}}(\XK,\XJ) + \FC_{\{k,i,j\}}(\XM,\XK,\XJ).
\end{aligned}
\end{equation}
Hence:
\begin{equation}\label{th403}
\begin{aligned}
\DKS \BGK_{\{j,k\}}(\XJ,\XM) &= \DJS \BGJ_{\{i,k\}}(\XK,\XM) = \DMS \BGM_{\{i,j\}}(\XK,\XJ),
\\
&\qquad \forall i ,j, k \quad {\rm such \; that} \quad 1 \leq i < j < k \leq N. 
\end{aligned}
\end{equation}
Next, proceed as before fixing values different from the "ground" values for
$\XK$, $\XJ$ and $\XM$ to obtain:
\begin{equation}\label{th404}  
\BGK_{\{j,k\}}(\XJ,\XM) = \chi_{i,j,k}  \DJS \DMS.
\end{equation} 
and:
\begin{equation}\label{th405}
\begin{aligned}
\DKS \DJS \DMS \chi_{i,j,k} = \DJS \DKS \DMS \chi_{j,i,k} = \cdots,
\end{aligned}
\end{equation}
so that:
\begin{equation}\label{th406}
\begin{aligned}
\chi_{i,j,k} =  \chi_{{\rm perm}(i,j,k)},
\end{aligned}
\end{equation}
where ${\rm perm}(i,j,k)$ denotes an arbitrary permutation of the subindexes $i$, $j$, and $k$.
Then, from \eqref{th402} the third order cliques may be obtained as:
\begin{equation}\label{th407}
\begin{aligned}
\FMS_{\{i,j,k\}}(\XK,\XJ) &= \; \chi_{i,j,k} \; \DKS \DJS \DMS + \FC_{\{i,j,k\}}(\XK,\XJ,\XM)
\end{aligned}
\end{equation}
\noindent
Hence, all the third order cliques for the mixed states distribution were obtained.
Proceed analogously to obtain all the higher order cliques.

\subsubsection{Joint distribution}\label{joint}

From the previous results the joint distribution may be obtained from the first, second and higher
order cliques as obtained from \eqref{th203}, \eqref{th306}, \eqref{th407} and its extensions,
as:
\begin{equation}\label{th501}
\begin{aligned}
\LOG \frac{\PMSX}{\PMSR} &= \; \sum_{k} \left[ \AK \DKS + \FCCK \right]
\\
&+ \; \sum_{<i,j>} \left[ \beta_{i,j} \; \DKS \DJS + \FC_{\{i,j\}}(\XK,\XJ) \right]
\\
&+ \; \sum_{<i,j,k>} \left[ \chi_{i,j,k} \; \DKS \DJS \DMS + \FC_{\{i,j,k\}}(\XK,\XJ,\XM) \right]
\\
&+ \cdots,
\end{aligned}
\end{equation}
\noindent
where $$\sum_{<i,j>} \equiv \sum_{\{i, j\} \in 2^{S}} = 
\stackblw{\sum_{i} \sum_{j}}{ 1 \le i < j \le N} = \sum_{i=1}^{N-1}\sum_{j=i+1}^{N},$$ 
and $$\sum_{<i,j,k>} \equiv \sum_{\{i, j, k \} \in 2^{S}} = 
\stackblw{\sum_{i} \sum_{j} \sum_{k}}{ 1 \le i < j < k \le N } = \sum_{i=1}^{N-2}\sum_{j=i+1}^{N-1}\sum_{k=j+1}^{N}. $$

Define:
\begin{equation}\label{th502}
\begin{aligned}
\VDX \equiv \FDX &= \; \sum_{i}  \AK \DKS 
+ \; \sum_{<i,j>}  \beta_{i,j} \; \DKS \DJS 
\\
&+ \; \sum_{<i,j,k>}  \chi_{i,j,k} \; \DKS \DJS \DMS 
+ \cdots.
\end{aligned}
\end{equation}
Hence, from \eqref{th3}, \eqref{th502}, and \eqref{th501} it results:
\begin{equation}\label{th503}
\begin{aligned}
\LOG \frac{\PMSX}{\PMSR} \equiv \VMX = \VDX \; + \; \VCX.
\end{aligned}
\end{equation}
To further stress the presence of the parameters, call $\ALPHAD$ 
a vector containing all the parameters of $\VDX$, see equation \eqref{th502}.
Similarly call $\ALPHAC$ the 
vector containing all the parameters of $\VCX$. Then, \eqref{th503} takes the form:
\begin{equation}\label{th504}
\begin{aligned}
\LOG \frac{\PMSX}{\PMSR} \; \equiv \; \FMSXP \; = \; \FDXP \; + \;  \FCXP.
\end{aligned}
\end{equation}
Define the partition functions:
\begin{equation}\label{th505}
\begin{aligned}
Z^{{a}}(\ALPHAC) \; = \; \int \; {\rm e}^{ \FCXP} \; m(d \X),
\end{aligned}
\end{equation}
\begin{equation}\label{th506}
\begin{aligned}
Z^{{d}}(\ALPHAD) \;  = \; \int \; {\rm e}^{ \FDXP}  \; m(d \X),
\end{aligned}
\end{equation}
and:
\begin{equation}\label{th507}
\begin{aligned}
\ZMS \; = \; \int {\rm e}^{\FMSXP}  \; m(d \X) \; = \; \int {\rm e}^{\FCXP +\FDXP} \; m(d \X)
\end{aligned}
\end{equation}
Define the generalized probability density function $\PDXP$ as:
\begin{equation}\label{th508}
\begin{aligned}
\PDXP = \frac{ {\rm e}^{\VDX}}{\ZD},
\end{aligned}
\end{equation}
so that:
\begin{equation}\label{th509}
\begin{aligned}
\ZD = \frac{ {1}}{p^{d}(\R ; \ALPHAD)}.
\end{aligned}
\end{equation}
Similarly, from \eqref{th3}:
\begin{equation}\label{th510}
\begin{aligned}
\ZC = \frac{ {1}}{p^{a}(\R ; \ALPHAC)}.
\end{aligned}
\end{equation}
Also, from \eqref{th504} and \eqref{th507}:
\begin{equation}\label{th511}
\begin{aligned}
\ZMS \;  = \; \frac{ {1}}{p^{m}(\R ; \ALPHAC)}.
\end{aligned}
\end{equation}
Collecting all these results in equation \eqref{th504}, it results:
\begin{equation}\label{th512}
\begin{aligned}
\PMSX \; = \; \frac{\PDXP \; \; \PCXP}{\ZN},
\end{aligned}
\end{equation}
where:
\begin{equation}\label{th513}
\begin{aligned}
\ZN \; = \; \frac{ \ZMS }{\ZD \; \; \ZC} = \int \PDXP \;\; \PCXP \;\; m(d \X),
\end{aligned}
\end{equation}
is the normalization constant, depending only on the parameters, for the product 
This product should not be misinterpreted as any independence condition.

From \eqref{th12}, \eqref{th11}, \eqref{th14}, \eqref{th303-1}, \eqref{th404} and its extensions:
\begin{equation}\label{th601}
\begin{aligned}
&\LFQP = \AK \; + \; \stackblw{\sum_{j}}{j \neq i} \beta_{i,j} \;  \DJS \; +
\\
& \qquad + \; \stackblw{\sum_{j,k}}{j \ne i, \; k \neq i, \; j < k} \chi_{i,j,k} \;  \DJS  \; \DMS
\; + \; \cdots \; \equiv \; \HKC.
\end{aligned}
\end{equation}
Hence:
\begin{equation}\label{th650}
\begin{aligned}
\QK \; = \; \frac{1}{1 +  \displaystyle \frac{{\rm e}^{\HKC}}{\PCRK}}.
\end{aligned}
\end{equation}
\\
\\
\subsection{The converse}\label{theconverse}

Conversely: if a MRF has the form given by equations \eqref{th512} and \eqref{th513}, satisfying theorem \ref{th:HC-1},
where $\PCXP$ is the pdf of an arbitrary absolutely continuous distribution function with respect
to Lebesgue measure, that is $\PCXP$ is a MRF by itself, i. e. satisfying equations \eqref{th3}, \eqref{th4} and theorem \ref{th:HC-1}, 
and,  $\PDXP$ is an arbitrary discrete random field given by
equations \eqref{th508}, \eqref{th502}, and \eqref{th509}, 
then, the conditional ms-pdf's have the form given by \eqref{th1} with the
$\QK$'s given by \eqref{th650}, where the $\HKC$'s are obtained from \eqref{th601}
using \eqref{th502}.
\\
Proof. Assume the joint generalized probability density function with respect to the measure $m$ given by 
equation \eqref{MEASURE} is given by:
\begin{equation}\label{conv001A}
\begin{aligned}
&\MPDFXV(\X) \; = \; \frac{{\rm e}^{\displaystyle \VMX}}{Z^{m}} = \frac{{\rm e}^{\displaystyle \VCX \; + \; \VDX}}{Z^{m}} =
\\
&= \frac{1}{Z^{m}} \; {\rm exp} \left( \CC + \; \sum_{i}  \AK \DKS + \; \sum_{<i,j>}  \beta_{i,j} \; \DKS \DJS \right.
\\
& \qquad \qquad + \left. \; \sum_{<i,j,k>}  \chi_{i,j,k} \; \DKS \DJS \DMS + \; \ldots \right).
\end{aligned}
\end{equation}
Note that $\MPDFXV(\X)  \neq 0$, for all $\X \in \mathcal{D}$.
\\
\noindent
Let $l \in S$, since $\MPDFXV(\X) \equiv p^{m}_{\XLRV, \XLCRV}(\XL,\XLC)$, then:
\begin{equation}\label{conv001B}
\begin{aligned}
\log \; &\frac{p^{m}_{\XLRV \; \vert \; \XLCRV}(\XL \; \vert \; \XLC)}{p^{m}_{\XLRV \; \vert \; \XLCRV}(\RL \; \vert \; \XLC)}
=
\log \; \frac{p^{m}_{\XLRV, \XLCRV}(\XL,\XLC)}{p^{m}_{\XLRV, \XLCRV}(\RL,\XLC)}
=
\\
&= \CCL + \; \UNOG^{*}_{\RL}(\XL) \; \Big(
 \alpha_{l} \; + 
\\
&\quad + \; \stackblw{\sum_{j}}{j \neq l} \beta_{l,j} \;  \UNOG^{*}_{\RJ}(\XJ)
\; +  \stackblw{\sum_{j,k}}{j \neq l, \; k \neq l, \; j < k} \chi_{l,j,k} \;  \UNOG^{*}_{\RJ}(\XJ)  \; \UNOG^{*}_{\RM}(\XM)
\; + \ldots \Big) =
\\
&= \; \UNOG^{*}_{\RL}(\XL) \; \; \HLC \; + \; \CCL,
\end{aligned}
\end{equation}
since $Q^{a}_{A} (\X_{A}) = 0$, for $\XL = \RL$ if $l \in A$, and where equation \eqref{th601} was used.

Hence:
\begin{equation}\label{CONV-101}
\begin{aligned}
p^{m}_{\XLRV \; \vert \; \XLCRV}&(\XL \; \vert \; \XLC)
= 
p^{m}_{\XLRV \; \vert \; \XLCRV}(\RL \; \vert \; \XLC)
\; 
{\rm e}^{\CCL} \; {\rm e}^{\displaystyle \UNOG^{*}_{\RL}(\XL) \; \; \HLC}
\\
&=
p^{m}_{\XLRV \; \vert \; \XLCRV}(\RL \; \vert \; \XLC)
\; 
{\rm e}^{\CCL} \; \Big( \UNOG_{\RL}(\XL) + \UNOG^{*}_{\RL}(\XL) \; {\rm e}^{\displaystyle \HLC} \Big)
\\
&=
p^{m}_{\XLRV \; \vert \; \XLCRV}(\RL \; \vert \; \XLC)
\; \Big( \UNOG_{\RL}(\XL) \; + \; \UNOG^{*}_{\RL}(\XL) \; {\rm e}^{\displaystyle \HLC} \; {\rm e}^{\CCL} \Big),
\end{aligned}
\end{equation}
because $\CCL = 0$, if $\XL = \RL$. Substituting equation \eqref{th4} in \eqref{CONV-101} it results:
\begin{equation}\label{CONV-111}
\begin{aligned}
p^{m}_{\XLRV \; \vert \; \XLCRV}(\XL \; \vert \; \XLC)
&= 
p^{m}_{\XLRV \; \vert \; \XLCRV}(\RL \; \vert \; \XLC)
\; \; \Big( \UNOG_{\RL}(\XL) \; + \;
\\
&\qquad + \; \UNOG^{*}_{\RL}(\XL) \; \frac{{\rm e}^{\displaystyle \HLC}}{\PCRL} \; \PCL \Big).
\end{aligned}
\end{equation}
Since $\int p^{m}_{\XLRV \; \vert \; \XLCRV}(\XL \; \vert \; \XLC) \; d m_{l} = 1$, and
\begin{equation}
\begin{aligned}
\int p^{m}_{\XLRV \; \vert \; \XLCRV}(\XL \; \vert \; \XLC) \; d m_{l}
&=  
\int p^{m}_{\XLRV \; \vert \; \XLCRV}(\XL \; \vert \; \XLC) \; d \MDX_{l} +
\int p^{m}_{\XLRV \; \vert \; \XLCRV}(\XL \; \vert \; \XLC) \; d \XL
\\
&=
p^{m}_{\XLRV \; \vert \; \XLCRV}(\RL \; \vert \; \XLC) + \int p^{m}_{\XLRV \; \vert \; \XLCRV}(\XL \; \vert \; \XLC) \; d \XL,
\end{aligned}
\end{equation}
then:
\begin{equation}\label{CONV-1001}
\begin{aligned}
1
= 
p^{m}_{\XLRV \; \vert \; \XLCRV}(\RL \; \vert \; \XLC) + \int p^{m}_{\XLRV \; \vert \; \XLCRV}(\XL \; \vert \; \XLC) \; d \XL.
\end{aligned}
\end{equation}
From equation \eqref{CONV-111}:
\begin{equation}
\begin{aligned}
\int &\frac{p^{m}_{\XLRV \; \vert \; \XLCRV}(\XL \; \vert \; \XLC)}{p^{m}_{\XLRV \; \vert \; \XLCRV}(\RL \; \vert \; \XLC)} \; d \X_{l}
\; = 
\\
&=
\int \UNOG_{\RL}(\XL) \; d \X_{l} \; + 
\; \int
{\NLFCL} \; \UNOG^{*}_{\RL}(\XL) \; {\rm e}^{\displaystyle \HLC}
\; d \X_{l}
\\
&=
\; 
\int
\UNOG_{\RL}(\XL) \; d \XL
+
{\rm e}^{\displaystyle \HLC}
\int 
\UNOG^{*}_{\RL}(\XL) \; {\NLFCL} \; d \XL.
\end{aligned}
\end{equation}
The first integral is zero since the integration set is a single point having Lebesgue measure zero, while the second integral is
a standard Lebesgue integral over the whole space except for a point of Lebesgue measure zero, so that:
\begin{equation}\label{CONV-1002}
\begin{aligned}
\int p^{m}_{\XLRV \; \vert \; \XLCRV}(\XL \; \vert \; \XLC) \; d \X_{l}
=
p^{m}_{\XLRV \; \vert \; \XLCRV}(\RL \; \vert \; \XLC)
\; 
{\rm e}^{\displaystyle \HLC}
\int 
\NLFCL \; d \XL.
\end{aligned}
\end{equation}
But:
\begin{equation}\label{CONV-1003}
\begin{aligned}
\int 
{\PCL} \; d \XL
= 1,
\end{aligned}
\end{equation}
so that, from equations \eqref{CONV-1001}, \eqref{CONV-1002}, and \eqref{CONV-1003}:
\begin{equation}\label{CONV-72}
\begin{aligned}
1 \; = \; p^{m}_{\XLRV \; \vert \; \XLCRV}(\RL \; \vert \; \XLC) \;
\Big(
1 \; + \; 
\frac{{\rm e}^{\displaystyle \HLC}}{p^{a}_{\XLRV \; \vert \; \XLCRV}(\RL \; \vert \; \XLC)} \; 
\Big).
\end{aligned}
\end{equation}
Call $\rho_l(\XLC) \; \equiv \; p^{m}_{\XLRV \; \vert \; \XLCRV}(\RL \; \vert \; \XLC)$, and $\rho_l^{*}(\XLC) = 1 - \rho_l(\XLC)$. 
Then, from equation \eqref{CONV-72}:
\begin{equation}
\begin{aligned}\label{CONV-2000}
\rho_l(\XLC)  
\; = \;
p^{m}_{\XLRV \; \vert \; \XLCRV}(\RL \; \vert \; \XLC)
\; = \;
\frac{1}
{
1 \; + \; 
\frac{\displaystyle {\rm e}^{\HLC}}{p^{a}_{\XLRV \; \vert \; \XLCRV}(\RL \; \vert \; \XLC)} \; 
}.
\end{aligned}
\end{equation}
Using equation \eqref{CONV-2000} in \eqref{CONV-111} one obtains:
\begin{equation}\label{CONV-777}
\begin{aligned}
p^{m}_{\XLRV \; \vert \; \XLCRV}(\XL \; \vert \; \XLC)
&= 
\rho_l(\XLC)
\; \UNOG_{\RL}(\XL) \; + \; \left(1 \; - \; \rho_l(\XLC) \right) \; \UNOG^{*}_{\RL}(\XL) \; \PCL
\\
&= 
\rho_l(\XLC)
\; \UNOG_{\RL}(\XL) \; + \; \rho_l^{*}(\XLC) \; \UNOG^{*}_{\RL}(\XL) \; \PCL,
\end{aligned}
\end{equation}
which is the desired result.

\section{Conclusions and Further Work} \label{conclusions}

A theoretical formulation of the mixed states random variable was presented, as well as a 
theoretical analysis for mixed states Markov Random Fields  
with probability mass  concentrated in a real value.
From the results obtained here, 
mainly equation \eqref{conv001A} which gives the joint ms-pdf with respect to the measure $m$,
and equations \eqref{CONV-2000} and \eqref{CONV-777} which give the conditional ms-pdf's
with respect to the measures $m_i$, previous results given in 
\cite{Bouthemy-1}, \cite{SPIE-06}, \cite{ICIP-06}, \cite{Hardouin-1}, \cite{Hardouin-2}, \cite{Gwenaelle-1}, \cite{Gwenaelle-2}, 
\cite{Gwenaelle-3}, \cite{Gwenaelle-4}, \cite{Gwenaelle-5}, \cite{yao-1}, are immediately obtained.

Equation \eqref{conv001A} permits to use the power of the Gibbs formulation using potentials to design the MRF, as an alternative to the use
of conditional distributions as was done in \cite{Bouthemy-1}, \cite{SPIE-06}, \cite{ICIP-06}, \cite{Hardouin-1}, \cite{Hardouin-2}, \cite{Gwenaelle-1}, \cite{Gwenaelle-2}, 
\cite{Gwenaelle-3}, \cite{Gwenaelle-4}, \cite{Gwenaelle-5}, \cite{yao-1}.

Results presented here will be extended in two directions.

In a sequel, the extension of the results given here
to Markov Random Fields of Mixed States variables which are mixtures of a denumerable set of probability 
mass concentrated in either label values and/or
multidimensional real values, and a standard absolutely continuous distributed multidimensional real random variable will be given.

The second direction corresponds to the analysis of Mixed States Markov Random Fields whose potentials 
present interaction between the "discrete" and the "absolutely continuous" distributed components of the mixtures for the potentials of the
Gibbs formulation, which are presently under study.

\break

\break

\tableofcontents

\end{document}